\documentclass[12pt]{article}%
\usepackage{amsfonts}
\usepackage{amsmath}
\usepackage{amssymb}
\usepackage{graphicx}%
\providecommand{\U}[1]{\protect\rule{.1in}{.1in}}
\newtheorem{theorem}{Theorem}

\newtheorem{lemma}[theorem]{Lemma}

\newtheorem{remark}[theorem]{Remark}

\newenvironment{proof}[1][Proof]{\noindent\textbf{#1.} }{\ \rule{0.5em}{0.5em}}
\begin{document}

\author{Bruce Anderson
\and J. Marshall Ash\thanks{Partially supported by NSF Grant DMS-9707011}
\and Roger Jones\thanks{Partially supported by NSF Grant DMS-9531526\noindent}
\and Daniel G. Rider
\and Bahman Saffari}
\title{Exponential sums with coefficients $0$ or $1$ and concentrated $L^{p}$ norms }
\maketitle

\begin{abstract}
A sum of exponentials of the form $f(x)=\exp\left(  2\pi iN_{1}x\right)
+\exp\left(  2\pi iN_{2}x\right)  +\cdot\cdot\cdot+\exp\left(  2\pi
iN_{m}x\right)  $, where the $N_{k}$ are distinct integers is called an
\textit{idempotent trigonometric polynomial} (because the convolution of $f$
with itself is $f$) or, simply, an \textit{idempotent}. We show that for every
$p>1,\,$and every set $E$ of the torus $\mathbb{T=R}/\mathbb{Z}$ with
$\left\vert E\right\vert >0,$ there are idempotents concentrated on $E$ in the
$L^{p}$ sense. More precisely, for each $p>1,$ there is an \textit{explicitly
calculated} constant $C_{p}>0$ so that for each $E$ with $\left\vert
E\right\vert >0$ and $\epsilon>0$ one can find an idempotent $f$ such that the
ratio $\left(  \int_{E}\left\vert f\right\vert ^{p}\left/  \int_{\mathbb{T}%
}\left\vert f\right\vert ^{p}\right.  \right)  ^{1/p}$ is greater than
$C_{p}-\epsilon$. This is in fact a lower bound result and, though
\textit{not} optimal, it is close to the best that \textit{our} method gives.
We also give both heuristic and computational evidence for the still open
problem of whether the $L^{p}$ concentration phenomenon fails to occur when
$p=1.$

\textbf{titre francais: }Sommes d'exponentielles \`{a} coefficients $0$ ou $1$
et concentration de normes $L^{p}$

\textbf{r\'{e}sum\'{e} francais: }Une somme d'exponentielles de la forme
$f(x)=\exp\left(  2\pi iN_{1}x\right)  +\exp\left(  2\pi iN_{2}x\right)
+\cdot\cdot\cdot+\exp\left(  2\pi iN_{m}x\right)  $, o\`{u} les $N_{k}$ sont
des entiers distincts, est appel\'{e}e un \textit{polyn\^{o}me
trigonom\'{e}trique idempotent } (car $f\ast f=f$) ou, simplement, un
\textit{idempotent}. Nous prouvons que pour tout r\'{e}el $p>1$, et tout
$E\subset\mathbb{T=R}/\mathbb{Z}$ avec $\left\vert E\right\vert >0,$ il existe
des idempotents concentr\'{e}s sur $E$ au sens de la norme $L^{p}$. Plus
pr\'{e}cis\'{e}ment, pour tout $p>1,$ nous \textit{calculons}
\textit{explicitement }une constante $C_{p}>0$ telle que pour tout $E$ avec
$\left\vert E\right\vert >0$, et tout r\'{e}el $\epsilon>0$, on puisse
construire un idempotent $f$ tel que le quotient $\left(  \int_{E}\left\vert
f\right\vert ^{p}\left/  \int_{\mathbb{T}}\left\vert f\right\vert ^{p}\right.
\right)  ^{1/p}$ soit sup\'{e}rieur \`{a} $C_{p}-\epsilon$. Ceci est en fait
un th\'{e}or\`{e}me de minoration qui, bien que \textit{non} optimal, est
proche du meilleur r\'{e}sultat que \textit{notre} m\'{e}thode puisse fournir.
Nous pr\'{e}sentons \'{e}galement des consid\'{e}rations heuristiques et aussi
num\'{e}riques concernant le probl\`{e}me (toujours ouvert) de savoir si le
ph\'{e}nom\`{e}ne de concentration $L^{p}$ a lieu ou non pour $p=1$.

\textbf{mots clefs: }idempotents, polyn\^{o}mes trigonom\'{e}triques
idempotents, normes $L^{p}$, noyau de Dirichlet, concentration de normes,
sommes d'exponentielles, conjecture de concentration en norme $L^{1}$,
op\'{e}rateurs faiblement restreints.

\textbf{keywords:} idempotents, idempotent trigonometric polynomials, $L^{p}$
norms, Dirichlet kernel, concentrating norms, sums of exponentials, $L^{1}$
concentration conjecture, weak restricted operators.

\textbf{classification code:} Primary 42A05; Secondary 42A10, 42A32.

\textbf{running title:} Exponential sums with coefficients 0 or 1

\end{abstract}

\bigskip

\section{Introduction\medskip}

\subsection{Concentrated $L^{p}$ norms}

Let $e\left(  x\right)  :=\exp\left(  2\pi ix\right)  .$ A sum of exponentials
of the form
\[
f\left(  x\right)  =\sum_{k=1}^{m}e\left(  N_{k}x\right)  \qquad\left(
x\in\mathbb{R}\right)  ,
\]
where the $N_{k}$ are distinct integers is called an \textit{idempotent
trigonometric polynomial }(because the convolution of $f$ with itself is $f$)
or, simply, an \textit{idempotent}. In the sequel we adopt the term
\textquotedblleft idempotent\textquotedblright\ for brevity, and we denote by
$\wp$ the set of all such idempotents:
\[
\wp:=\{\sum_{n\in S}e\left(  nx\right)  :S\text{ is a finite set of
non-negative integers}\}.
\]
The simplest example of an $f\in\wp$ is (one form of) the \textit{Dirichlet
Kernel} of length $n,$ defined by
\begin{equation}
D_{n}(x):=\sum_{\nu=0}^{n-1}e(\nu x)=\frac{\sin\left(  n\pi x\right)  }%
{\sin\left(  \pi x\right)  }\cdot e^{i\left(  n-1\right)  \pi x}. \label{1}%
\end{equation}
Consider any function $g\in L^{p}\left(  \mathbb{T}\right)  ,$ where
$\mathbb{T=R}/\mathbb{Z}$, and any set $E\subset\mathbb{T}$ with $\left\vert
E\right\vert >0.$ (Throughout, $\left\vert E\right\vert $ denotes the Lebesgue
measure of $E$.) If
\begin{equation}
\left(  \int_{E}|g\left(  x\right)  |^{p}dx\left/  \int_{\mathbb{T}}|g\left(
x\right)  |^{p}dx\right.  \right)  ^{1/p}\geq\alpha\label{2}%
\end{equation}
(where $1\leq p<\infty$ and $0<\alpha<1),$ we say that \textquotedblleft%
\textit{at least a proportion }$\alpha$\textit{\ of the }$L^{p}$\textit{\ norm
of }$g $\textit{\ is concentrated on }$E"$ or, equivalently, that
\textquotedblleft\textit{the function }$g$\textit{\ has }$L^{p}$%
\textit{\ concentration }$\geq\alpha$\textit{\ on }$E$\textit{.}%
\textquotedblright\ We will now recall a challenging (and still partially
open) problem on idempotents which can be expressed in terms of this notion of
$L^{p}$ concentration.\medskip

About 25 years ago it was discovered that for any (\textit{arbitrarily small})
arc $J$ of the torus $\mathbb{T}$ with $\left\vert J\right\vert >0,$ there
always exists an \textit{idempotent }$f$ with \textit{at least }%
$48\%$\textit{\ of its }$L^{2}$\textit{\ norm concentrated on }$J$. The origin
of this curious fact occurred about 1977 when one of us (J. M. Ash) was
attempting to show that an operator $T$ defined on $L^{2}(\mathbb{T})$,
commuting with translations, and of restricted weak type $(2,2)$ is
necessarily a bounded operator on $L^{2}\left(  \mathbb{T}\right)  .$ That $T$
is of restricted weak type $(2,2)$ means that there is a constant $C=C\left(
T\right)  >0$ such that, for every characteristic function $\chi$ of a subset
of $\mathbb{T},$
\[
\sup_{\alpha>0}\left(  \alpha^{2}\operatorname*{measure}\{x\in\mathbb{T}%
:|T\chi(x)|>\alpha\}\right)  \leq C\left\Vert \chi\right\Vert _{L^{2}}^{2}.
\]
Ash\cite{As} was only able to show that this condition was equivalent to there
being \textit{some} positive amount of $L^{2}$ concentration for idempotents.
More explicitly, define the absolute constant $C_{2}^{\ast}$ as the largest
real number such that for \textit{every} arc $J\subset\mathbb{T}$ with
$\left\vert J\right\vert >0$, we have the inequality
\begin{equation}
\sup_{f\in\wp}\left(  \int_{J}|f\left(  x\right)  |^{2}dx\left/
\int_{\mathbb{T}}|f\left(  x\right)  |^{2}dx\right.  \right)  ^{1/2}\geq
C_{2}^{\ast}. \label{3}%
\end{equation}
Thus the issue was whether $C_{2}^{\ast}$ was $0$ or positive. Luckily, at
just the same time, Michael Cowling\thinspace\cite{Co} proved, by another
method, that every commuting with translations weak restricted type $(2,2)$
operator\textit{\ is} necessarily a bounded operator on $L^{2}.$ (Actually
Cowling\thinspace\cite{Co} proved more. His result allowed the underlying
group to be \textit{any} amenable group, not just $\mathbb{T}$.) This proved,
of course, that $C_{2}^{\ast}$ was indeed positive, but did not give any
effective estimate for it. However, a series of concrete estimates quickly
followed. The referee of \cite{As} obtained $C_{2}^{\ast}\geq.01$, S.
Pichorides\thinspace\cite{Pi} obtained $C_{2}^{\ast}\geq.14$, H. L.
Montgomery\thinspace\cite{Mo} and J.-P. Kahane\thinspace\cite{Ka2} obtained
several better lower bounds. (The ideas of H. L. Montgomery were
\textquotedblleft deterministic\textquotedblright\ while those of J.-P. Kahane
used probabilistic methods from \cite{Ka1}.) Finally, in \cite{AJS}, three of
us achieved the lower bound
\begin{equation}
\gamma_{2}:=\max_{x>0}\frac{\sin x}{\sqrt{\pi x}}=.4802..., \label{4}%
\end{equation}
which, in \cite{DPQ1}, was proved to be best possible. (See \cite{DPQ2} for a
more detailed exposition of the contents of \cite{DPQ1}.)\medskip

To get a little more feel for what to expect, let $\zeta$ be any point of
density (also called \textquotedblleft Lebesgue point\textquotedblright) of a
set $E\subset\mathbb{T}.$ Then for every $p\in\lbrack1,\infty\lbrack$ the
sequence of functions $\{g_{n}\},$ where $g_{n}(x):=D_{n}(x-\zeta)=\sum
_{\nu=0}^{n-1}e(-\nu\zeta)\cdot e(\nu x),$ have $L^{p}$ concentration tending
to $1$ as $n\rightarrow\infty.$ However, the trigonometric polynomials $g_{n}$
are not idempotents, since the non-zero coefficients are not all equal to $1.$
Note, however, that all the coefficients do have \textit{modulus} $1.$ The
difficulty of the matters studied in \cite{AJS} and \cite{DPQ1}, as well as
those of the present paper lies precisely in the fact that the trigonometric
polynomials $f\in\wp$ have all their coefficients equal to $0$ or $1,$ which
is a very drastic constraint.\medskip

At this stage we make an obvious remark: in all the $L^{2}$ concentration
problems on (small) arcs of $\mathbb{T}$ studied in \cite{AJS} and \cite{DPQ1}
and in all their $L^{p}$ analogues studied in the present paper, it is
\textit{equivalent} to work on arcs of $\mathbb{T}$ or on intervals of
$\left[  0,1\right]  .$ We usually find it convenient to use \textquotedblleft
arcs of $\mathbb{T}$\textquotedblright\ in statements of theorems, but
\textquotedblleft intervals of $\left[  0,1\right]  $\textquotedblright\ in
their proofs!

\subsection{The $L^{2}$ and $L^{p}$ problems}

The results in \cite{AJS} and \cite{DPQ1} were satisfying but, as usual, they
led to further questions. The first two were:

\begin{description}
\item[a)] Can we replace \textquotedblleft arc\textquotedblright\ (or
\textquotedblleft interval\textquotedblright) with \textquotedblleft set of
positive measure?\textquotedblright

\item[b)] Can we replace $L^{2}$ with $L^{p}$ for any $p\geq1$?
\end{description}

For each $p\in\lbrack1,\infty\lbrack,$ define $C_{p}$ as the largest number
such that for every set $E$, $E\subset\mathbb{T}$ with $\left\vert
E\right\vert >0$, the inequality
\begin{equation}
\sup_{f\in\wp}\left\Vert f\right\Vert _{L^{p},E}\left/  \left\Vert
f\right\Vert _{L^{p}}\right.  :=\sup_{f\in\wp}\left(  \int_{E}|f|^{p}dx\left/
\int_{\mathbb{T}}|f|^{p}dx\right.  \right)  ^{1/p}\geq C_{p} \label{*5}%
\end{equation}
holds. Similarly, define $C_{p}^{\ast}$ as the largest number such that for
every arc $J$, $J\subset\mathbb{T}$, the inequality
\begin{equation}
\sup_{f\in\wp}\left\Vert f\right\Vert _{L^{p},J}\left/  \left\Vert
f\right\Vert _{L^{p}}\right.  :=\sup_{f\in\wp}\left(  \int_{J}|f|^{p}dx\left/
\int_{\mathbb{T}}|f|^{p}dx\right.  \right)  ^{1/p}\geq C_{p}^{\ast} \label{*6}%
\end{equation}
holds. The definitions of $C_{p}$ and $C_{p}^{\ast}$ are extended to the limit
cases $p=\infty$ in the usual way. Obviously,%
\begin{equation}
C_{p}\leq C_{p}^{\ast}. \label{*7}%
\end{equation}

With regard to question \textbf{a)}, although the definitions allow the
possibility for $C_{2}$ to be smaller than $C_{2}^{\ast}$, in \cite{AJS} it is
shown that both are equal to the constant $\gamma_{2}$ defined in (\ref{4}).
Whatever the value of $p\geq1$, there is no result in this paper which changes
when the supremum is taken over all sets of positive measure rather than over
all arcs. So we \textit{conjecture} that inequality (\ref{*7}) is in fact an
equality:%
\begin{equation}
C_{p}=C_{p}^{\ast}, \label{*8}%
\end{equation}
although we have no proof of this except for $p=2$ and $p=\infty$.

Question \textbf{b)} is harder. In \cite{DPQ1}, the constant $\gamma_{2}$
defined in (\ref{4}) is shown to be a lower bound for every $C_{p}$ when
$p\geq2$. This is not altogether satisfying for two reasons. First, the cases
$1\leq p<2$ are not addressed. Second, since the constant function $1$ is in
$\wp$, and $\left\Vert 1\right\Vert _{L^{\infty},A}\left/  \left\Vert
1\right\Vert _{L^{\infty}}\right.  =1$ for \textit{any} non-empty set
$A\subset\mathbb{T},$ so that $C_{\infty}=1$, one might hope to show that
$\lim_{p\rightarrow\infty}C_{p}=1.$ In section \ref{1.3} we state new results
for the $L^{p}$ cases, together with some remaining open problems.

\subsection{\label{1.3}Statement of the result}

This paper is devoted to proving one single theorem, Theorem \ref{one} below.
It was announced in the Comptes Rendus note \cite{AAJRS} (in a weaker form,
and presented as two distinct results). The aim of this paper is to supply the
proofs of \cite{AAJRS}, to strengthen the first result thereof, and to unify
the results in the form of a single theorem (which is valid for all sets of
positive measure and not just for all arcs). Our theorem is stated in terms of
the \textquotedblleft constants\textquotedblright\ $c_{p}$ and $c_{p}^{\ast},$
defined as follows:%
\[
c_{p}:=\sup_{0<\omega<1/2}\frac{\sin\left(  \pi\omega\right)  /\left(
\pi\omega\right)  }{2^{1+1/p}\left(  \left\lfloor 1/\omega\right\rfloor
+1+\frac{1}{p-1}\left(  \frac{3}{8}\right)  ^{p}\left\lfloor 1/\omega
\right\rfloor \right)  ^{1/p}},
\]
where for a real number $r,$
\begin{align*}
\left\lceil r\right\rceil  &  =\text{ceiling of }r=\text{the smallest integer
greater than or equal to }r\text{,}\\
\left\lfloor r\right\rfloor  &  =\text{floor of }r=\text{the greatest integer
less than or equal to }r;
\end{align*}
and%
\[
c_{p}^{\ast}:=\left(  \frac{2}{\pi^{p+1}}\int_{0}^{\infty}\left\vert
\frac{\sin x}{x}\right\vert ^{p}dx\right)  ^{1/p}\cdot\max_{0\leq\omega\leq
1}\frac{\sin\left(  \pi\omega\right)  }{\omega^{1-1/p}}.
\]
As $p$ increases from $1$ to $+\infty$ (resp. from $2$ to $+\infty$), $c_{p}$
(resp. $c_{p}^{\ast}$) increases from $0$ to $.5$ (resp. from $\gamma_{2}$
$=0.48\dots$ to $1$). That $c_{p}^{\ast}$ tends to $1$ as $p\rightarrow\infty$
follows from an easy calculation, which is done in Remark \ref{r:1} for the
reader's convenience.\medskip

\begin{theorem}
\label{one}Whenever $1<p<\infty,$ we have the estimate%
\[
C_{p}\geq\left\{
\begin{array}
[c]{lll}%
c_{p} & \  & \text{ if }1<p\leq2\\
&  & \\
c_{p}^{\ast} &  & \text{ if }2\leq p<\infty
\end{array}
\right.  .
\]
In other words, if $p>1$ and $\epsilon>0$ are given, then for each set
$E\subset\mathbb{T}$ with $\left\vert E\right\vert >0$, there is a finite set
of integers $S=S(E,p,\epsilon)$ such that
\begin{equation}
\int_{E}\left\vert \sum_{n\in S}e(nx)\right\vert ^{p}dx\left/
%TCIMACRO{\dint _{\mathbb{T}}}%
%BeginExpansion
{\displaystyle\int_{\mathbb{T}}}
%EndExpansion
\left\vert \sum_{n\in S}e(nx)\right\vert ^{p}dx\right.  \geq\left\{
\begin{array}
[c]{lll}%
c_{p}^{p}-\epsilon & \  & \text{ if }1<p\leq2\\
&  & \\
c_{p}^{\ast p}-\epsilon &  & \text{ if }2\leq p<\infty
\end{array}
\right.  . \label{interval}%
\end{equation}
Furthermore, $c_{p}^{\ast}$ (and a fortiori $C_{p}$) tends to $1$ as $p$ tends
to infinity.
\end{theorem}

For a slightly larger lower estimate of $C_{p},$ see inequality
(\ref{giratio1}) at the end of Section 2 (where Part I of Theorem \ref{one} is
proved).\medskip

\begin{remark}
\label{r:1}The estimate in Part II of Theorem \ref{one}, which is quite good
although not optimal, does have two virtues. First, it is sharp when $p=2,$
since
\[
\frac{2}{\pi}\int_{0}^{\infty}\left\vert \frac{\sin x}{x}\right\vert
^{2}dx=1.
\]
Second, it implies that $\lim_{p\rightarrow\infty}C_{p}=1$. Indeed
\begin{align*}
\liminf_{p\rightarrow\infty}C_{p}  &  \geq\liminf_{p\rightarrow\infty}%
c_{p}^{\ast}\\
&  =\liminf_{p\rightarrow\infty}\left(  \max_{0\leq\omega\leq1}\frac{\sin
\pi\omega}{\pi\omega^{1-1/p}}\right)  \lim_{p\rightarrow\infty}\left(
\frac{2}{\pi}\right)  ^{1/p}\lim_{p\rightarrow\infty}\left(  \int_{0}^{\infty
}\left\vert \frac{\sin x}{x}\right\vert ^{p}dx\right)  ^{1/p}\\
&  \geq\lim_{p\rightarrow\infty}\frac{\sin\pi\left(  1/p\right)  }{\pi\left(
1/p\right)  ^{1-1/p}}\cdot1\cdot\operatorname*{Essup}\left\vert \frac{\sin
x}{x}\right\vert =\lim_{p\rightarrow\infty}\frac{\sin\left(  \pi/p\right)
}{\left(  \pi/p\right)  }\lim_{p\rightarrow\infty}\left(  1/p\right)
^{1/p}=1.
\end{align*}

\end{remark}

To prove Theorem \ref{one}, it obviously suffices to prove the inequalities
$C_{p}\geq c_{p}$ (for all $p>1$) and $C_{p}\geq c_{p}^{\ast}$ (for all
$p\geq2$). These two inequalities will be proved, respectively, in the next
sections 2 and 3.\medskip

As for the \textit{open problems} pertaining to the case $p=1$ (conjectures of
non-concentration in the $L^{1}$ sense), we shall state them in Section 4 at
the end of the paper.

\section{Proof of Theorem 1; Part I: $C_{p}\geq c_{p}$ (for all $p>1$)}

\noindent\textbf{2.1.} Since the proof is quite technical and\ computational,
before giving the full proof we start by sketching a (heuristic) outline of
the beginning of the proof.\medskip

\noindent\textbf{Outline of (the beginning of) the proof. \ }Let \thinspace$q$
be a large odd positive integer and let $m$ be an integer at least as large.
We begin with the special case of $J=[\frac{1}{q}-\frac{1}{mq},\frac{1}%
{q}+\frac{1}{mq}],$ where $m\geq q$. The idea of the proof is to think of
$\mathbb{T}$ (or rather of some suitable interval of length $1,$ which is more
convenient in the proofs) as \textquotedblleft partitioned\textquotedblright%
\ (except for common endpoints) into $q$ congruent arcs of the form $\left[
\frac{2\nu-1}{2q},\frac{2\nu+1}{2q}\right]  ,$ with centers at $\frac{\nu}%
{q},\left(  \nu=0,1,\dots,q-1\right)  $ and common length $\frac{1}{q}.$ The
idempotent $D_{mq}(qx),$ where $D_{n}\left(  x\right)  $ is the Dirichlet
kernel defined by the relation (\ref{1}) of the introduction, has period $1/q$
and when restricted to $\left[  -\frac{1}{2}\frac{1}{q},\frac{1}{2}\frac{1}%
{q}\right]  $ behaves approximately like the Dirac measure. Now consider the
idempotent $D_{\omega q}(x),$ where $\omega\in]0,1/2[$ is chosen to make
$\omega q$ an integer. This, when restricted to a small neighborhood of the
set $\left\{  \frac{0}{q},\frac{1}{q},\frac{2}{q},...,\frac{q-1}{q}\right\}
,$ behaves roughly (as far as its \textit{modulus} is concerned) like the
function $1/x.$ Thus the idempotent $\varpi(x)=$ $D_{mq}(qx)D_{\omega q}(x)$
satisfies $\int_{0}^{1}|\varpi(x)|^{p}dx\approx\sum_{j=1}^{q}1/j^{p}$ and
$\int_{J}|\varpi(x)|^{p}dx\approx1^{-p}.$ Concentration at $1/q$ follows since
$\sum_{j=1}^{q}j^{-p}$ is bounded.\medskip

This outline will be made rigorous in the proof below. We will also have to
treat the problem of concentration at points which are \textit{not} of the
form $1/q.\medskip$

\noindent\textbf{2.2. }To prove Theorem \ref{one} we need the following lemma,
which will be used in the proof of Part II as well.

\begin{lemma}
\label{l:1}\textit{Let} $D_{N}(x)$\textit{\ be the Dirichlet Kernel defined
by} (\ref{1}), so that $\left\vert D_{N}\left(  x\right)  \right\vert
=\left\vert \sin\left(  \pi Nx\right)  /\sin\left(  \pi x\right)  \right\vert
,$ \textit{and let} $p$\textit{\ be greater than }$1.$ \textit{Then}
\[
\int_{0}^{1}\left\vert D_{N}(x)\right\vert ^{p}dx=\delta_{p}N^{p-1}%
+o_{p}(N^{p-1}),\quad\quad N\rightarrow\infty
\]
\textit{where}
\[
\delta_{p}=\frac{2}{\pi}\int_{0}^{\infty}\left\vert \frac{\sin u}%
{u}\right\vert ^{p}du
\]
\textit{and} $o_{p}$ \textit{is the \textquotedblleft little
o\textquotedblright\ notation of Landau modified to emphasize the dependence
of the associated constant on }$p.$
\end{lemma}

\textit{More precisely,}%
\[
\int_{0}^{1}|D_{N}(u)|^{p}du=\delta_{p}N^{p-1}+R_{p}\left(  N\right)  ,
\]
\textit{where the error term} $R_{p}\left(  N\right)  $ \textit{satisfies:}
\[
R_{p}\left(  N\right)  =\left\{
\begin{array}
[c]{lll}%
O_{p}\left(  N^{p-3}\right)  &  & \text{ if }p>3\\
O\left(  \log N\right)  &  & \text{ if }p=3\\
O_{p}\left(  1\right)  &  & \text{ if }1<p<3
\end{array}
.\right.
\]
%

%TCIMACRO{\TeXButton{Proof}{\proof}}%
%BeginExpansion
\proof
%EndExpansion
This result is classical, but we give the full proof for the reader's
convenience. Since $D_{N}$ is even, we need only estimate $2\int_{0}%
^{1/2}\left\vert D_{N}(x)\right\vert ^{p}dx.$ By the triangle inequality, this
differs from
\[
A_{p}\left(  N\right)  :=2\int_{0}^{1/2}|\frac{\sin N\pi u}{\pi u}|^{p}du
\]
by at most
\[
E_{p}\left(  N\right)  :=2\int_{0}^{1/2}|\sin N\pi u|^{p}\left(  \frac{1}%
{\sin^{p}\pi u}-\frac{1}{\left(  \pi u\right)  ^{p}}\right)  du.
\]

Substituting $x:=N\pi u$ yields
\begin{align*}
A_{p}\left(  N\right)   &  =\frac{2}{\pi N}\int_{0}^{N\pi/2}|\frac{\sin
x}{x/N}|^{p}dx\\
&  =\left(  \frac{2}{\pi}\int_{0}^{\infty}|\frac{\sin x}{x}|^{p}dx\right)
N^{p-1}-\frac{2}{\pi}N^{p-1}\int_{N\pi/2}^{\infty}|\frac{\sin x}{x}|^{p}dx\\
&  =\delta_{p}N^{p-1}+O_{p}\left(  1\right)  ,
\end{align*}
since
\[
\int_{N\pi/2}^{\infty}|\frac{\sin x}{x}|^{p}dx<\int_{N\pi/2}^{\infty}%
x^{-p}dx=\frac{\pi^{-p+1}2^{p-1}}{p-1}N^{-\left(  p-1\right)  }.
\]
So proving the lemma reduces to proving that%
\[
E_{p}\left(  N\right)  =\left\{
\begin{array}
[c]{lll}%
O_{p}\left(  N^{p-3}\right)  &  & \text{ if }p>3\\
O\left(  \log N\right)  &  & \text{ if }p=3\\
O_{p}\left(  1\right)  &  & \text{ if }1<p<3
\end{array}
\right.  .
\]
For $0<u\leq1/2$, the inequality%
\begin{align*}
\frac{1}{\sin^{p}\left(  \pi u\right)  }-\frac{1}{\left(  \pi u\right)  ^{p}}
&  \leq\frac{\left(  \pi u\right)  ^{p}-\sin^{p}\left(  \pi u\right)
}{\left(  \pi u\right)  ^{p}\sin^{p}\left(  \pi u\right)  }\\
&  =O_{p}\left(  1\right)  \frac{u^{p}\left(  1-\left(  1+O\left(
u^{2}\right)  \right)  ^{p}\right)  }{u^{2p}}=O_{p}\left(  u^{2-p}\right)
\end{align*}
immediately leads to the two estimates%
\begin{equation}
|\sin N\pi u|^{p}\left(  \frac{1}{\sin^{p}\pi u}-\frac{1}{\left(  \pi
u\right)  ^{p}}\right)  =\left(  N\pi u\right)  ^{p}O_{p}\left(
u^{2-p}\right)  =O_{p}\left(  N^{p}u^{2}\right)  \label{est1}%
\end{equation}
and%

\begin{equation}
|\sin N\pi u|^{p}\left(  \frac{1}{\sin^{p}\pi u}-\frac{1}{\left(  \pi
u\right)  ^{p}}\right)  =1\cdot O_{p}\left(  u^{2-p}\right)  =O_{p}\left(
u^{2-p}\right)  . \label{est3}%
\end{equation}
If $1<p<3,$ from estimate (\ref{est3}) we have%
\[
E_{p}\left(  N\right)  \leq\int_{0}^{1/2}O_{p}\left(  u^{2-p}\right)
du=O_{p}\left(  1\right)  ;
\]
if $p=3,$ from estimates (\ref{est1}) and (\ref{est3}) we have%
\begin{align*}
E_{3}\left(  N\right)   &  \leq\int_{0}^{1/N}O\left(  N^{3}u^{2}\right)
du+\int_{1/N}^{1/2}O\left(  u^{2-3}\right)  du\\
&  \leq O\left(  N^{3}N^{-3}\right)  +O\left(  \log N\right)  =O\left(  \log
N\right)  ;
\end{align*}
and if $p>3,$ again from estimates (\ref{est1}) and (\ref{est3}) we have%
\begin{align*}
E_{p}\left(  N\right)   &  \leq\int_{0}^{1/N}O_{p}\left(  N^{p}u^{2}\right)
du+\int_{1/N}^{1/2}O_{p}\left(  u^{2-p}\right)  du\\
&  =O_{p}\left(  N^{p}N^{-3}\right)  +O_{p}\left(  N^{-\left(  3-p\right)
}\right)  =O_{p}\left(  N^{p-3}\right)  .
\end{align*}
Thus the lemma is proved.\medskip

\noindent\textbf{2.3. }We now proceed to prove Theorem \ref{one} in detail. We
find it convenient to split the proof into eight steps. The first three steps
deal with concentration at $1/q$ and the remaining five steps with the general
case.\medskip

\noindent\underline{\textbf{First Step.}}\textbf{\ Concentration at }%
$1/q$\textbf{: lower estimation of numerator.}

Let $\omega=\omega(q)$ be a constant in $(0,1/2)$ such that $\omega q$ is an
integer. We will estimate
\begin{equation}
\left(  \int_{J}|\varpi(x)|^{p}dx\left/  \int_{0}^{1}|\varpi(x)|^{p}dx\right.
\right)  ^{1/p}, \label{ratio}%
\end{equation}
where $\varpi(x):=D_{mq}(qx)D_{\omega q}(x).\medskip$

We begin with the numerator, $N,$ of (\ref{ratio}).\medskip

Suppose that $|\frac{1}{q}-x|\leq\frac{1}{mq},$ and let $\delta:=\frac{1}%
{q}-x.$ Then, since $\sin x$ has a bounded derivative,
\[
|D_{\omega q}(x)|=|D_{\omega q}(\frac{1}{q}-\delta)|=\left|  \frac{\sin
\pi\omega+O\left(  \delta q\right)  }{\sin\frac{\pi}{q}+O\left(
\delta\right)  }\right|  =\left|  \frac{\sin\pi\omega}{\sin\frac{\pi}{q}%
}\right|  +O\left(  1\right)
\]
\[
=\left|  \frac{\sin\pi\omega}{\frac{\pi}{q}}\right|  +O\left(  1\right)  ,
\]
since $m\geq q.$ Using this and Minkowski's inequality in the form
\newline$\left\|  F+G\right\|  _{L^{p}}\geq\left\|  F\right\|  _{L^{p}%
}-\left\|  G\right\|  _{L^{p}},$ we have
\begin{equation}
N\geq\left(  \int_{J}|D_{mq}(qx)|^{p}\left(  \frac{q\sin\pi\omega}{\pi
}\right)  ^{p}dx\right)  ^{1/p}-O\left(  \left(  \int_{J}|D_{mq}%
(qx)|^{p}dx\right)  ^{1/p}\right)  . \label{num}%
\end{equation}
Substitute $u=qx$ to get
\begin{equation}
N\geq\left\{  \left(  \frac{\sin\pi\omega}{\pi}\right)  q^{1-1/p}-O\left(
q^{-1/p}\right)  \right\}  \left(  \int_{-1/m}^{1/m}|D_{mq}(u)|^{p}du\right)
^{1/p}. \label{num1}%
\end{equation}
Define $\Delta:=\left(  \int_{0}^{1}|D_{mq}(u)|^{p}du\right)  ^{1/p}.$ Since
$D_{mq}$ is even, we may write
\[
\int_{-1/m}^{1/m}|D_{mq}(u)|^{p}du=\Delta^{p}-2\int_{1/m}^{1/2}|D_{mq}%
(u)|^{p}du.
\]
Use the estimates $\left|  \sin mq\pi u\right|  \leq1$ and $\left|  \sin\pi
u\right|  \geq2u$ to control the last integral, thereby obtaining the
estimate
\[
\int_{-1/m}^{1/m}|D_{mq}(u)|^{p}du=\Delta^{p}-O\left(  m^{p-1}\right)  .
\]
By the lemma, $\Delta^{p}\simeq\delta_{p}q^{p-1}m^{p-1},$ which together with
(\ref{num1}) implies
\begin{equation}
N\geq\Delta\left\{  \left(  \frac{\sin\pi\omega}{\pi}\right)  q^{1-1/p}%
-O\left(  q^{-1/p}\right)  \right\}  \left\{  1-O\left(  q^{-1+1/p}\right)
\right\}  , \label{num2}%
\end{equation}
or, more simply,
\begin{equation}
N\geq\Delta\left(  \frac{\sin\pi\omega}{\pi}\right)  q^{1-1/p}\left\{
1-o(1)\right\}  . \label{num3}%
\end{equation}

\noindent\underline{\textbf{Second Step.}}\textbf{\ Concentration at }%
$1/q$\textbf{: upper estimation of denominator.\ }

Passing now to the estimate of the denominator, $D,$ of (\ref{ratio}), we
have
\begin{equation}
D^{p}=\int_{0}^{1}|D_{mq}(qx)|^{p}\left|  D_{\omega q}(x)\right|  ^{p}dx.
\label{denom1}%
\end{equation}
We now estimate this in great detail. We decompose
\[
D^{p}=\sum_{j=0}^{q-1}\int_{\frac{j}{q}}^{\frac{j+1}{q}}\left|  D_{mq}%
(qx)\right|  ^{p}\left|  D_{\omega q}(x)\right|  ^{p}dx.
\]
Let $x=y+\frac{j}{q},dy=dx,$ to get
\[
D^{p}=\sum_{j=0}^{q-1}\int_{0}^{\frac{1}{q}}\left|  D_{mq}(qy+j)\right|
^{p}\left|  D_{\omega q}\left(  y+\frac{j}{q}\right)  \right|  ^{p}dy.
\]
Since $\left|  D_{mq}(qy+j)\right|  =\left|  D_{mq}(qy)\right|  $ for any
integer $j,$%
\[
D^{p}=\sum_{j=0}^{q-1}\int_{0}^{\frac{1}{q}}\left|  D_{mq}(qy)\right|
^{p}\left|  D_{\omega q}\left(  y+\frac{j}{q}\right)  \right|  ^{p}dy.
\]
Let $t=qy$ to get
\[
D^{p}=\frac{1}{q}\sum_{j=0}^{q-1}\int_{0}^{1}\left|  D_{mq}(t)\right|
^{p}\left|  D_{\omega q}\left(  \frac{t}{q}+\frac{j}{q}\right)  \right|
^{p}dt.
\]
Interchange sum and integral:
\begin{equation}
D^{p}=\int_{0}^{1}\left|  D_{mq}(t)\right|  ^{p}\frac{1}{q}\sum_{j=0}%
^{q-1}\left|  D_{\omega q}\left(  \frac{t}{q}+\frac{j}{q}\right)  \right|
^{p}dt. \label{denom2}%
\end{equation}
Replace the sum by its supremum over all $t\in\lbrack0,1]$ which is the same
as
\[
\sup_{x\in\lbrack0,\frac{1}{q})}\sum_{j=0}^{q-1}\left|  D_{\omega q}\left(
x+\frac{j}{q}\right)  \right|  ^{p},
\]
so recalling that $\Delta^{p}=\int_{0}^{1}\left|  D_{mq}(t)\right|  ^{p}dt,$
we have
\[
D^{p}\leq\Delta^{p}\frac{1}{q}\sup_{x\in\lbrack0,\frac{1}{q})}\sum_{j=0}%
^{q-1}\left|  D_{\omega q}\left(  x+\frac{j}{q}\right)  \right|  ^{p}.
\]
Since $D_{\omega q}$ is even and has period $1,$%
\[
D^{p}\leq2\Delta^{p}\frac{1}{q}\sup_{x\in\lbrack0,\frac{1}{q})}\sum
_{j=0}^{\frac{q-1}{2}}\left|  D_{\omega q}\left(  x+\frac{j}{q}\right)
\right|  ^{p}.
\]
Break the sum into two pieces using the standard estimates $\left|  D_{\omega
q}(x+j/q)\right|  \leq\omega q$ when $j\leq\left\lfloor 1/\omega\right\rfloor
$ and
\[
\sup_{x\in\lbrack0,\frac{1}{q})}\left|  D_{\omega q}(x+\frac{j}{q})\right|
\leq1/\sin\left(  \frac{\pi j}{q}\right)  \leq\frac{q}{2j}
\]
when $j\in\left[  \left\lfloor 1/\omega\right\rfloor +1,(q-1)/2\right]  .$ We
have
\begin{align*}
D^{p}  &  \leq2\Delta^{p}\frac{1}{q}\left\{  \sum_{j=0}^{\left\lfloor \frac
{1}{\omega}\right\rfloor }\left(  \omega q\right)  ^{p}+\sum_{j=\left\lfloor
\frac{1}{\omega}\right\rfloor +1}^{\frac{q-1}{2}}\left(  \frac{q}{2j}\right)
^{p}\right\} \\
\  &  \leq2\Delta^{p}\frac{1}{q}\left\{  \left(  \left\lfloor \frac{1}{\omega
}\right\rfloor +1\right)  \left(  \omega q\right)  ^{p}+\left(  \frac{q}%
{2}\right)  ^{p}\int_{\left\lfloor \frac{1}{\omega}\right\rfloor }^{\infty
}\frac{dx}{x^{p}}\right\} \\
\  &  =2\Delta^{p}q^{p-1}\omega^{p}\left\{  \left(  \left\lfloor \frac
{1}{\omega}\right\rfloor +1\right)  +\frac{\left\lfloor \frac{1}{\omega
}\right\rfloor }{p-1}\left(  \frac{1}{2}\right)  ^{p}\left(  \frac{\frac
{1}{\omega}}{\left\lfloor \frac{1}{\omega}\right\rfloor }\right)  ^{p}\right\}
\\
&  =2\Delta^{p}q^{p-1}\omega^{p}\left\{  \left\lfloor \frac{1}{\omega
}\right\rfloor +1+\frac{\left\lfloor \frac{1}{\omega}\right\rfloor \rho}%
{p-1}\right\}  ,
\end{align*}
where
\[
\rho=\left(  \frac{1}{2}\frac{1/\omega}{\left\lfloor 1/\omega\right\rfloor
}\right)  ^{p}.
\]

\noindent\underline{\textbf{Third Step.}}\textbf{\ Concentration at }%
$1/q$\textbf{: conclusion.}

Now combine this estimate with estimate (\ref{num3}) to get
\[
\frac{N}{D}\geq\frac{\Delta\left(  \frac{\sin\pi\omega}{\pi}\right)
q^{1-1/p}\left\{  1-o(1)\right\}  }{\left(  2\Delta^{p}q^{p-1}\omega
^{p}\left\{  \left\lfloor \frac{1}{\omega}\right\rfloor +1+\frac{\left\lfloor
\frac{1}{\omega}\right\rfloor \rho}{p-1}\right\}  \right)  ^{1/p}}%
\]
or
\[
\frac{N}{D}\geq\frac{\left(  \frac{\sin\pi\omega}{\pi\omega}\right)  \left\{
1-o(1)\right\}  }{\left(  2\left\{  \left\lfloor \frac{1}{\omega}\right\rfloor
+1+\frac{\left\lfloor \frac{1}{\omega}\right\rfloor \rho}{p-1}\right\}
\right)  ^{1/p}}.
\]
The numbers $\omega=\omega\left(  q\right)  $ appearing in the above two steps
are, by construction, rational numbers. However it is clear, from their
definition, that as $q$ varies these $\omega\left(  q\right)  $ are everywhere
dense in $\left[  0,1/2\right]  .$ (Cf. also eighth step of this proof.) As
will be made explicit below, this estimate, when extended to general
intervals, is sufficient to prove Theorem 1. (Notice that $\rho<\left(
3/8\right)  ^{p}$ since $\omega<1/2$.)\medskip

\noindent\underline{\textbf{Fourth Step.}} \textbf{General case: heuristic
search for the concentrated exponential sum.}

From now on, our goal is to extend the above estimate (of the third step) to
any interval. This fourth step is purely heuristic.\medskip

Given any interval $J\subset]0,1[,$ we can find $q$ (and in fact infinitely
many such $q$'s) so that for some integer $k\in\lbrack1,q[,$ $[\frac{k}%
{q}-\frac{1}{qm},\frac{k}{q}+\frac{1}{qm}]\subset J,$ where $m:=q$. (This
choice of $m$ is for technical reasons that will become clear in the eighth
step of the proof.) It is convenient to pick $q$ \textit{prime,} which implies
that $k$ and $q$ are relatively prime. Hence there exists a unique pair
$\left(  a,b\right)  $\ of integers such that
\begin{equation}
ak-bq=1,\quad\left(  0<a<q\text{ \quad and \quad}0<b<k\right)  \label{inverse}%
\end{equation}
i.e., (the conjugate class of) $a$ is the multiplicative inverse of $k$ in the
finite field $Z_{q}=GF\left(  q\right)  .$ Multiplication by $a,$ when reduced
modulo $q,$ defines a bijection $\alpha$ from $\left\{  0,1,...,q-1\right\}  $
to itself. Furthermore, $\alpha(k)=1.$ For the sake of this heuristic
argument, temporarily suppose that (as previously) $\omega$ is chosen so that
$\omega q$ is an integer. (In fact, in the rigorous argument below, we shall
choose \ $\omega$ according to another criterion, so that $\omega q$ will
\textit{not }be an integer, but instead of $\omega q$ we shall use the integer
$\left\lceil \omega q\right\rceil .$ Now the idempotent $D_{\omega q}(ax)$
behaves very much like the idempotent $D_{\omega q}(x),$ except that the
former does at $k/q$ what the latter does at $1/q$. Since $D_{mq}(qx)$ is
constant and large on the set $\left\{  \frac{0}{q},\frac{1}{q},\frac{2}%
{q},...,\frac{q-1}{q}\right\}  ,$ and $D_{\omega q}(ax)$ restricted to this
set takes on the same set of values as $D_{\omega q}(x)$ did, but has its
maximum at $k/q$ (instead of at $1/q);$ it seems reasonable that, the
idempotent $D_{mq}(qx)D_{\omega q}(ax)$ should work here. The definition of
the idempotent $G\left(  x\right)  $\ analyzed below was motivated by these
considerations.\medskip

\noindent\underline{\textbf{Fifth Step.}}\ \textbf{General case: the desired
concentrated exponential sum.}

Let $E$ be a subset of $\mathbb{T}$ of positive measure and let $\epsilon>0$
be given. First we will find an integer $Q$, and an $\eta=\eta\left(
Q,\epsilon\right)  >0;$ then an interval $J$ of the form $J=\left[  \frac
{k}{q}-\frac{1}{mq},\frac{k}{q}+\frac{1}{mq}\right]  $ so that $\left\vert
J\cap E\right\vert >\left(  1-\eta\right)  \left\vert J\right\vert $; and
finally we will define an idempotent $f$ depending on $k$ and $Q$ which is
$\epsilon$-close to being sufficiently concentrated on first $J$ and then $E$.
Actually in what follows we will always take $m$ to be equal to $q$, but we
leave $m$ in the calculations since some increase of the concentration
constants may be available by taking other values of $m$.

First we define $\eta$. The function $f$ will have the form%
\[
f\left(  x\right)  =G\left(  x\right)  \sum_{n=0}^{mQ-1}e\left(  nqx\right)
\text{.}
\]
Since $G$ will turn out to be a sum of $\leq q$ exponentials (see
(\ref{defG}), $\left\Vert G\right\Vert _{\infty}\leq q$ ; so the decomposition
$J=\left(  J\smallsetminus E\right)  \cup\left(  J\cap E\right)  $ allows
\begin{align*}
\int_{J\setminus E}\left\vert f\right\vert ^{p}dx  &  \leq\eta\left\vert
J\right\vert \left(  \left\Vert f\right\Vert _{\infty}\right)  ^{p}\\
&  \leq\eta\frac{2}{mq}\left(  qmQ\right)  ^{p}\\
&  =\left(  \eta Q\right)  2q^{p-1}\left(  mQ\right)  ^{p-1}.
\end{align*}
In the sixth step below, we will get
\[
\int_{J}\left\vert f\right\vert ^{p}dx\geq\left(  \frac{\sin\pi\omega}{\pi
}\right)  ^{p}q^{p-1}\delta_{p}\left(  mQ\right)  ^{p-1}+o\left(
q^{p-1}\left(  mQ\right)  ^{p-1}\right)
\]
as long as $q$ is large enough. Pick $\eta\left(  Q,\epsilon\right)  $ so
small that from this will follow
\begin{equation}
\int_{J\cap E}\left\vert f\right\vert ^{p}dx=\int_{J}\left\vert f\right\vert
^{p}dx-\int_{J\setminus E}\left\vert f\right\vert ^{p}dx\geq\left(
1-\epsilon\right)  ^{p}\left(  \frac{\sin\pi\omega}{\pi}\right)  ^{p}%
q^{p-1}\delta_{p}\left(  mQ\right)  ^{p-1}. \label{Efix1}%
\end{equation}

Now that we know how to choose $\eta$, we show how to find $J=J\left(
E,\eta\right)  $. We will use $m=q$ while choosing $J$. Almost every point
$\xi$ has the property that there are infinitely many primes $q$ and integers
$k$ for which
\begin{equation}
\left\vert \xi-\frac{k}{q}\right\vert <\frac{1}{q^{2}}\text{. (See \cite{H}.)}
\label{close}%
\end{equation}
We may assume that $E$ has an irrational point of density $\xi$ for which
inequality (\ref{close}) holds for infinitely many primes $q$. Suppose that
the prime $q$ is so large that $K=\left[  \xi-\frac{2}{q^{2}},\xi+\frac
{2}{q^{2}}\right]  $ satisfies $\left\vert K\diagdown E\right\vert \leq
\frac{\eta}{2}\left\vert K\right\vert $ and such that condition (\ref{close})
holds. With $J=\left[  \frac{k}{q}-\frac{1}{q^{2}},\frac{k}{q}+\frac{1}{q^{2}%
}\right]  $, we have
\begin{align*}
\left\vert J\diagdown E\right\vert  &  \leq\left\vert K\diagdown E\right\vert
\leq\\
&  \leq\frac{\eta}{2}\left\vert K\right\vert =\frac{\eta}{2}\frac{4}{q^{2}%
}=\eta\frac{2}{q^{2}}=\eta\left\vert J\right\vert ,
\end{align*}
so that $\left\vert J\cap E\right\vert >\left(  1-\eta\right)  \left\vert
J\right\vert $, as required.

Let $a$ be (uniquely) defined by (\ref{inverse}), which in turn uniquely
defines the bijection from $\left\{  0,1,...,q-1\right\}  $ into itself
(reduction modulo $q$ of multiplication by $a$). We have $\alpha\left(
r\right)  =ra-sq$ where $s$ is the largest (necessarily non-negative) integer
such that $sq\leq ra.$ This leads us to consider the following sets
$E_{j},j=1,2,...$. For each integer $j\geq0$ let $E_{j}$ denote the set of
those $r\in\left\{  0,1,...,q-1\right\}  $ such that $\alpha\left(  r\right)
=ra-jq.$ \textit{A priori} the $E_{j}$ are pairwise disjoint, and it is
straightforward to check that
\[
E_{j}=\left\{
\begin{array}
[c]{ccc}%
\left\{  r\in\mathbb{N}:t_{j}\leq r<t_{j+1}\right\}  &  & \text{ if }0\leq
j<a\\
&  & \\
\emptyset &  & \text{ if }j\geq a
\end{array}
\right.  ,
\]
where $\left\{  t_{j}\right\}  $ is the (\textit{strictly increasing}) finite
sequence of integers defined by $t_{0}:=0,$ $t_{j}=\left\lceil
jq/a\right\rceil $ if $0<j<a,$ and $t_{j}=q$ if $j=a.\medskip$

Thus $\left\{  E_{j}\right\}  _{0\leq j<a}$ is a partition of $\left\{
0,1,...,q-1\right\}  ,$ and we have $\alpha\left(  r\right)  =ra-jq$ when
$r\in E_{j}.\medskip$

Now pick $\omega\in\left]  0,1/2\right[  $ so that $\omega a$ is an integer
$\ell$ at our disposal (with the obvious constraint $0<\ell<a/2).$ Instead of
the \textquotedblleft heuristic\textquotedblright\ idempotent $D_{\omega
q}\left(  x\right)  $ suggested in the Fourth Step above, we now consider the
idempotent
\begin{equation}
G(x):=\sum_{r=0}^{t_{\ell}-1}e(\alpha(r)x), \label{defG}%
\end{equation}
where, in view of the above calculations, $t_{\ell}=\left\lceil \ell
q/a\right\rceil =\left\lceil \omega q\right\rceil .$ To make good use of this
$G\left(  x\right)  ,$ we now need to perform a long calculation:
\[
G(x)=\sum_{j=0}^{\ell-1}\sum_{r=t_{j}}^{t_{j+1}-1}e(\alpha(r)x)=\sum
_{j=0}^{\ell-1}\sum_{r=t_{j}}^{t_{j+1}-1}e((ra-jq)x)
\]%
\[
=\sum_{j=0}^{\ell-1}e(-jqx)\sum_{s=0}^{t_{j+1}-t_{j}-1}e((sa+t_{j}a)x)
\]%
\[
=\sum_{j=0}^{\ell-1}e(-jqx)e(t_{j}ax)\sum_{s=0}^{t_{j+1}-t_{j}-1}%
e(sax)=\sum_{j=0}^{\ell-1}e((t_{j}a-jq)x)\left(  \frac{e((t_{j+1}-t_{j}%
)ax)-1}{e(ax)-1}\right)
\]%
\[
=\frac{1}{e(ax)-1}\sum_{j=0}^{\ell-1}e(-jqx)\left(  e(t_{j+1}ax)-e(t_{j}%
ax)\right)
\]%
\[
=\frac{1}{e(ax)-1}\left(  \sum_{j=1}^{\ell}e(-(j-1)qx)e(t_{j}ax)-\sum
_{j=0}^{\ell-1}e(-jqx)e(t_{j}ax)\right)
\]%
\[
=\frac{1}{e(ax)-1}\left(  \sum_{j=1}^{\ell}e(-jqx)e(qx)e(t_{j}ax)-\sum
_{j=0}^{\ell-1}e(-jqx)e(t_{j}ax)\right)
\]%
\[
=\frac{1}{e(ax)-1}\left\{  \sum_{j=1}^{\ell-1}e(-jqx)e(t_{j}ax)\left(
e(qx)-1\right)  +e(-\ell qx)e(qx)e(t_{\ell}ax)-1\right\}
\]%
\[
=\frac{e(qx)-1}{e(ax)-1}\left\{  \sum_{j=1}^{\ell-1}e(-jqx)\left(
e(t_{j}ax)-1\right)  +\sum_{j=1}^{\ell-1}e(-jqx)\right\}
\]%
\[
+\frac{1}{e(ax)-1}\left\{  e(t_{\ell}ax)\left(  e(-\ell qx)e(qx)-1\right)
+e(t_{\ell}ax)-1\right\}
\]%
\[
=\frac{e(qx)-1}{e(ax)-1}\left\{  \sum_{j=1}^{\ell-1}e(-jqx)\left(
e(t_{j}ax)-1\right)  +e(-qx)\frac{e(-(\ell-1)qx)-1}{e(-qx)-1}\right\}
\]%
\[
+\frac{1}{e(ax)-1}\left\{  e(t_{\ell}ax)\left(  e(-(\ell-1)qx)-1\right)
+e(t_{\ell}ax)-1\right\}
\]%
\[
=\frac{e(qx)-1}{e(ax)-1}\left\{  \sum_{j=1}^{\ell-1}e(-jqx)\left(
e(t_{j}ax)-1\right)  \right\}
\]%
\[
+\frac{e(qx)-1}{e(ax)-1}e(-qx)\frac{e(-(\ell-1)qx)-1}{e(-qx)-1}+e(t_{\ell
}ax)\left(  \frac{e(-(\ell-1)qx)-1}{e(ax)-1}\right)  +\frac{e(t_{\ell}%
ax)-1}{e(ax)-1}
\]%
\[
=\frac{e(qx)-1}{e(ax)-1}\left\{  \sum_{j=1}^{\ell-1}e(-jqx)\left(
e(t_{j}ax)-1\right)  \right\}
\]%
\[
+\frac{e(-(\ell-1)qx)-1}{e(ax)-1}\left(  -1+e(t_{\ell}ax)\right)
+\frac{e(t_{\ell}ax)-1}{e(ax)-1}
\]%
\[
=\frac{e(qx)-1}{e(ax)-1}\left\{  \sum_{j=1}^{\ell-1}e(-jqx)\left(
e(t_{j}ax)-1\right)  \right\}  +\frac{e(t_{\ell}ax)-1}{e(ax)-1}\left(
e(-(\ell-1)qx)-1+1\right)  .
\]
Hence we may write
\[
G(x)=\frac{e(qx)-1}{e(ax)-1}\sum_{j=1}^{\ell-1}e(-jqx)\left(  e(t_{j}%
ax)-1\right)  +\frac{e(t_{\ell}ax)-1}{e(ax)-1}e(-(\ell-1)qx)=:G_{1}+G_{2}.
\]
Let $S:=\{qn+\alpha(r):r=0,1,...,t_{\ell}-1$ and $n=0,1,...,mq-1\}.$ Also let
\[
f(x):=\sum_{j\in S}e(jx)=\sum_{n=0}^{mQ-1}e(nqx)G(x).
\]
This will be the desired concentrated sum of exponentials. Thus for each
$p>1,$ we must estimate the ratio
\begin{equation}
r:=\left(  \int_{J\cap E}\left\vert f(x)\right\vert ^{p}dx\left/  \int_{0}%
^{1}\left\vert f(x)\right\vert ^{p}dx\right.  \right)  ^{1/p}. \label{girat}%
\end{equation}

\noindent\underline{\textbf{Sixth Step:}}\textbf{\ General case: lower
estimation of the numerator.}

We now estimate the numerator of the ratio (\ref{girat}).\medskip

As in the fourth step, we may assume without loss of generality that $J\ $is
centered at $\frac{k}{q}$ and has length $\frac{2}{mq}$, where $m=q$. If $x\in
J,$ write $x=:\frac{k}{q}+y,$ so that $\left\vert y\right\vert $ $\leq\frac
{1}{mq}.$ Then
\[
G(x)=\sum_{r=0}^{t_{\ell}-1}e(\frac{r}{q})e(\alpha(r)y),
\]
so that
\[
G(x)=\sum_{r=0}^{t_{\ell}-1}e(\frac{r}{q})+\sum_{r=0}^{t_{\ell}-1}e(\frac
{r}{q})O(\frac{1}{m})=\sum_{r=0}^{t_{\ell}-1}e(\frac{r}{q})+O(\frac{t_{\ell}%
}{m})=\frac{e(t_{\ell}/q)-1}{e(1/q)-1}+O(\frac{\omega q}{m}).
\]
Thus on $J$ we have
\begin{align*}
\left\vert G(x)\right\vert  &  =\left\vert \frac{\sin\pi\left(  t_{\ell
}/q\right)  }{\sin\pi\left(  1/q\right)  }\right\vert +O(\frac{\omega q}%
{m})=\left\vert \frac{\sin\pi\left(  \omega+\left(  \left\lceil \omega
q\right\rceil -\omega q\right)  /q\right)  }{\sin\pi\left(  1/q\right)
}\right\vert +O(\frac{\omega q}{m})\\
&  =\frac{\sin\pi\omega}{\pi}q+O(1)
\end{align*}
since $m\geq q$. Also $\sum_{j=0}^{mQ-1}e(qjx)=D_{mQ}(qx),$ so applying
Minkowski's inequality to the numerator in the ratio (\ref{girat}), we have
\begin{align*}
\left(  \int_{J}\left\vert f(x)\right\vert ^{p}dx\right)  ^{1/p}  &  =\left(
\int_{J}\left\vert \frac{\sin\pi\omega}{\pi}q+O(1)\right\vert ^{p}\left\vert
D_{mQ}(qx)\right\vert ^{p}dx\right)  ^{1/p}\\
&  \geq\frac{\sin\pi\omega}{\pi}q\left(  \int_{J}\left\vert D_{mQ}%
(qx)\right\vert ^{p}dx\right)  ^{1/p}-O\left(  \int_{J}\left\vert
D_{mQ}(qx)\right\vert ^{p}dx\right)  ^{1/p}.
\end{align*}
The same reasoning that led to equation (\ref{num2}) above now produces the
following estimate for the numerator of (\ref{girat}):
\[
\left(  \int_{J}\left\vert f(x)\right\vert ^{p}dx\right)  ^{1/p}\geq\frac
{\sin\pi\omega}{\pi}q^{1-1/p}\Delta-O(q^{-1/p}\Delta),
\]
where $\Delta=\left(  \int_{0}^{1}\left\vert D_{mQ}(u)\right\vert
^{p}du\right)  ^{1/p}$. Taking (\ref{Efix1}) into account brings us to
\begin{equation}
\left(  \int_{J\cap E}\left\vert f\right\vert ^{p}dx\right)  ^{1/p}\geq\left(
1-\epsilon\right)  \frac{\sin\pi\omega}{\pi}q^{1-1/p}\delta_{p}^{1/p}\left(
mQ\right)  ^{1-1/p}\text{.} \label{ginum}%
\end{equation}

\noindent\underline{\textbf{Seventh Step.}}\textbf{\ General case: upper
estimation of the denominator.}

We now estimate the denominator of the ratio (\ref{girat}). To do this we will
need the following Lemmas.

\begin{lemma}
\label{l:4}Let $p>1$, $0<\theta<1$, and%
\[
K_{\theta,N}=\left\{  \left(  x,y\right)  \in\mathbb{Z}^{2}:x+\theta^{-1}y\leq
N,x\geq0,y\geq0\right\}  .
\]
Then for arbitrary $N\geq4$,
\[
\int_{0}^{1}\int_{0}^{1}\left\vert \sum_{\left(  m,n\right)  \in K_{\theta,N}%
}e\left(  mx+ny\right)  \right\vert ^{p}dxdy\leq C_{p}N^{2p-2}
\]
uniformly with respect to $\theta$ and $N$.
\end{lemma}

See \cite{As1} for the proof of this. That proof extends (and is much indebted
to) work of Yudin and Yudin\cite{YY} for $p=1$ to higher values of $p$. (Cf.
\cite{NP}.)

\begin{lemma}
\label{l:y}Let $p>1$, $0<\theta<1$, and%
\[
L_{\theta,N}=\left\{  \left(  x,y\right)  \in\mathbb{Z}^{2}:x+\theta^{-1}y\leq
N,x>0,y\geq0\right\}  .
\]
Then for arbitrary $N\geq4$,
\[
\int_{0}^{1}\int_{0}^{1}\left\vert \sum_{\left(  m,n\right)  \in L_{\theta,N}%
}e\left(  mx+ny\right)  \right\vert ^{p}dxdy\leq C_{p}N^{2p-2}
\]
uniformly with respect to $\theta$ and $N$.
\end{lemma}

\begin{proof}
This is immediate from the last Lemma because $K_{\theta,N}$ is the disjoint
union of $L_{\theta,N}$ and $\left\{  \left(  0,n\right)  :n=0,1,\dots
,\left\lfloor \theta N\right\rfloor \right\}  $; and by Lemma \ref{l:1},
\[
\int_{0}^{1}\int_{0}^{1}\left\vert \sum_{n=0}^{\left\lfloor \theta
N\right\rfloor }e\left(  0x+ny\right)  \right\vert ^{p}dxdy=\int_{0}%
^{1}\left\vert D_{\left\lfloor \theta N\right\rfloor }\left(  y\right)
\right\vert ^{p}dy\ll N^{p-1}<N^{2p-2}.
\]

\end{proof}

To study the denominator of the ratio (\ref{girat}), because of Minkowski's
inequality
\begin{equation}
\left\Vert f\right\Vert _{p}\leq\left\Vert \frac{\sin\pi mqQx}{\sin\pi
qx}G_{1}\right\Vert _{p}+\left\Vert \frac{\sin\pi mqQx}{\sin\pi qx}%
G_{2}\right\Vert _{p}=:A+B, \label{gidenom0}%
\end{equation}
it is enough to estimate these last two terms separately. For the first of
these we have
\[
A^{p}=\int_{0}^{1}\left\vert \frac{\sin\pi mqQx}{\sin\pi qx}\right\vert
^{p}\left\vert \frac{\sin\pi qx}{\sin\pi ax}\right\vert ^{p}\left\vert
\sum_{n=1}^{\ell-1}e(-nqx)\left(  e(t_{n}ax)-1\right)  \right\vert ^{p}dx
\]%
\begin{align*}
&  =\int_{0}^{1}\left\vert \sin\pi mqQx\right\vert ^{p}\left\vert \sum
_{n=1}^{\ell-1}e(-nqx)\frac{e(t_{n}ax)-1}{\sin\pi ax}\right\vert ^{p}dx\\
&  =\int_{0}^{1}\left\vert \sin\pi mqQx\right\vert ^{p}\left\vert \sum
_{n=1}^{\ell-1}e(-nqx)D_{t_{n}}\left(  ax\right)  \right\vert ^{p}dx\\
&  \leq\int_{0}^{1}\left\vert \sum_{n=1}^{\ell-1}e(-nqx)D_{t_{n}}\left(
ax\right)  \right\vert ^{p}dx=\sum_{j=0}^{q-1}\int_{j/q}^{\left(  j+1\right)
/q}\left\vert \sum_{n=1}^{\ell-1}e(-nqx)D_{t_{n}}\left(  ax\right)
\right\vert ^{p}dx\\
&  =\sum_{j=0}^{q-1}\int_{0}^{1/q}\left\vert \sum_{n=1}^{\ell-1}%
e(-nx-ajn)D_{t_{n}}\left(  ax+aj/q\right)  \right\vert ^{p}dx\\
&  =\sum_{j=0}^{q-1}\int_{0}^{1/q}\left\vert \sum_{n=1}^{\ell-1}%
e(-nx)D_{t_{n}}\left(  ax+aj/q\right)  \right\vert ^{p}dx\\
&  =\sum_{j=0}^{q-1}\int_{0}^{1/q}\left\vert \sum_{n=1}^{\ell-1}%
e(-nx)D_{t_{n}}\left(  ax+j/q\right)  \right\vert ^{p}dx\\
&  =\frac{1}{q}\sum_{j=0}^{q-1}\int_{0}^{1}\left\vert \sum_{n=1}^{\ell
-1}e(-ny/q)D_{t_{n}}\left(  ay/q+j/q\right)  \right\vert ^{p}dy
\end{align*}
where the last five steps are justified by the identity $\left[  0,1\right]
=\cup_{j=0}^{q-1}\left[  j/q,\left(  j+1\right)  /q\right]  $, the
substitutions $x\rightarrow x+j/q$, the periodicity of $e\left(  x\right)  $
being $1$, the fact that $t\rightarrow at$ is a one to one correspondence on
$\mathbb{Z}_{q}$, and the substitution $x=y/q$. Summarizing,%
\[
A^{p}\leq\frac{1}{q}\sum_{j=0}^{q-1}\int_{0}^{1}\left\vert P_{y}\left(
ay/q+j/q\right)  \right\vert ^{p}dy
\]
where $P_{y}\left(  x\right)  =\sum_{n=1}^{\ell-1}e(-ny)D_{t_{n}}\left(
x\right)  $. But a theorem of Marcinkiewicz and Zygmund \cite{Z} asserts%
\[
\frac{1}{q}\sum_{j=0}^{q-1}\left\vert P_{y}\left(  ay/q+j/q\right)
\right\vert ^{p}\ll\int_{0}^{1}\left\vert P_{y}\left(  x\right)  \right\vert
^{p}dx.
\]
We may now apply Lemma \ref{l:y} since the index set of $P_{y}\left(
x\right)  $ is $L_{\theta,N}$ with $\theta=a/q$ and $N=\left(  \ell-1\right)
q/a$ to get
\[
A^{p}\ll\int_{0}^{1}\int_{0}^{1}\left\vert P_{y}\left(  x\right)  \right\vert
^{p}dxdy\ll q^{2p-2},
\]%
\begin{equation}
A<Cq^{2-2/p}. \label{gidenom1}%
\end{equation}

Passing to the other term, we have
\[
B^{p}=\int_{0}^{1}\left\vert \frac{\sin\pi mqQx}{\sin\pi qx}\right\vert
^{p}\left\vert \frac{\sin\pi t_{\ell}ax}{\sin\pi ax}\right\vert ^{p}\left\vert
e(-(\ell-1)qx)\right\vert ^{p}dx
\]%
\[
=\int_{0}^{1}\left\vert D_{mQ}(qx)\right\vert ^{p}\left\vert D_{\left\lceil
\omega q\right\rceil }(ax)\right\vert ^{p}dx.
\]
This is estimated in a very similar way to the way $D$ was estimated above,
but with a couple of twists.
\[
B^{p}=\sum_{j=0}^{q-1}\int_{\frac{j}{q}}^{\frac{j+1}{q}}\left\vert
D_{mQ}(qx)\right\vert ^{p}\left\vert D_{\left\lceil \omega q\right\rceil
}(ax)\right\vert ^{p}dx.
\]
Let $x=y+\frac{j}{q},dy=dx,$ to get
\[
B^{p}=\sum_{j=0}^{q-1}\int_{0}^{\frac{1}{q}}\left\vert D_{mQ}(qy+j)\right\vert
^{p}\left\vert D_{\left\lceil \omega q\right\rceil }\left(  ay+\frac{ja}%
{q}\right)  \right\vert ^{p}dy.
\]
Since $\left\vert D_{mQ}(qy+j)\right\vert =\left\vert D_{mQ}(qy)\right\vert $
for any integer $j,$%
\[
B^{p}=\sum_{j=0}^{q-1}\int_{0}^{\frac{1}{q}}\left\vert D_{mQ}(qy)\right\vert
^{p}\left\vert D_{\left\lceil \omega q\right\rceil }\left(  ay+\frac{ja}%
{q}\right)  \right\vert ^{p}dy.
\]

Let $t=qy$ to get
\[
B^{p}=\frac{1}{q}\sum_{j=0}^{q-1}\int_{0}^{1}\left\vert D_{mQ}(t)\right\vert
^{p}\left\vert D_{\left\lceil \omega q\right\rceil }\left(  \frac{at}{q}%
+\frac{ja}{q}\right)  \right\vert ^{p}dt.
\]
Interchange sum and integral:
\[
B^{p}=\int_{0}^{1}\left\vert D_{mQ}(t)\right\vert ^{p}\frac{1}{q}\sum
_{j=0}^{q-1}\left\vert D_{\left\lceil \omega q\right\rceil }\left(  \frac
{at}{q}+\frac{ja}{q}\right)  \right\vert ^{p}dt.
\]
Now as $j$ varies between $0$ and $q-1,$ so does $ja,$ modulo $q.$ Thus the
inner sum of this term may be written as
\[
\sum_{j=0}^{q-1}\left\vert D_{\left\lceil \omega q\right\rceil }\left(
\frac{at}{q}+\frac{j}{q}\right)  \right\vert ^{p}.
\]
The supremum of this over all $t\in\lbrack0,1]$ is the same as
\[
\sup_{x\in\lbrack0,\frac{1}{q})}\sum_{j=0}^{q-1}\left\vert D_{\left\lceil
\omega q\right\rceil }\left(  x+\frac{j}{q}\right)  \right\vert ^{p},
\]
so recalling that $\Delta^{p}=\int_{0}^{1}\left\vert D_{mQ}(t)\right\vert
^{p}dt,$ we have
\[
B^{p}\leq\Delta^{p}\frac{1}{q}\sup_{x\in\lbrack0,\frac{1}{q})}\sum_{j=0}%
^{q-1}\left\vert D_{\left\lceil \omega q\right\rceil }\left(  x+\frac{j}%
{q}\right)  \right\vert ^{p}.
\]
Since $\left\vert D_{\left\lceil \omega q\right\rceil }\right\vert $ is even
and has period $1,$%
\[
B^{p}\leq2\Delta^{p}\frac{1}{q}\sup_{x\in\lbrack0,\frac{1}{q})}\sum
_{j=0}^{\frac{q-1}{2}}\left\vert D_{\left\lceil \omega q\right\rceil }\left(
x+\frac{j}{q}\right)  \right\vert ^{p}.
\]
Break the sum into two pieces using the standard estimates $\left\vert
D_{\left\lceil \omega q\right\rceil }(x)\right\vert \leq\left\lceil \omega
q\right\rceil $ when $j\leq\left\lfloor \frac{1}{\omega}\right\rfloor $
and$,$
\[
\sup_{x\in\lbrack0,\frac{1}{q})}\left\vert D_{\left\lceil \omega q\right\rceil
}(x+\frac{j}{q})\right\vert \leq1/\sin\left(  \frac{\pi j}{q}\right)
\leq\frac{q}{2j}
\]
when $j\in\left[  \left\lfloor 1/\omega\right\rfloor ,\left(  q-1\right)
/2\right]  .$ We have
\begin{align*}
B^{p}  &  \leq2\Delta^{p}\frac{1}{q}\left\{  \sum_{j=0}^{\left\lfloor \frac
{1}{\omega}\right\rfloor }\left\lceil \omega q\right\rceil ^{p}+\sum
_{j=\left\lfloor \frac{1}{\omega}\right\rfloor +1}^{\frac{q-1}{2}}\left(
\frac{q}{2j}\right)  ^{p}\right\} \\
\  &  \leq2\Delta^{p}\frac{1}{q}\left\{  \left(  \left\lfloor 1/\omega
\right\rfloor +1\right)  \left\lceil \omega q\right\rceil ^{p}+\left(
\frac{q}{2}\right)  ^{p}\int_{\left\lfloor 1/\omega\right\rfloor }^{\infty
}\frac{dx}{x^{p}}\right\} \\
\  &  =2\Delta^{p}q^{p-1}\omega^{p}\left\{  \left(  \left\lfloor
1/\omega\right\rfloor +1\right)  \left(  \left\lceil \omega q\right\rceil
/\omega q\right)  ^{p}+\frac{\left\lfloor 1/\omega\right\rfloor }{p-1}\left(
\frac{1}{2}\right)  ^{p}\left(  \frac{1/\omega}{\left\lfloor 1/\omega
\right\rfloor }\right)  ^{p}\right\}  .
\end{align*}
Clearly $\left\lceil \omega q\right\rceil /\left(  \omega q\right)  \leq2,$
whence
\begin{equation}
B^{p}\leq2^{p+1}\Delta^{p}q^{p-1}\omega^{p}\left(  \left\lfloor 1/\omega
\right\rfloor +1+\frac{\left\lfloor 1/\omega\right\rfloor }{p-1}\rho\right)  ,
\label{gidenom2}%
\end{equation}
where
\[
\rho=\left(  \frac{1}{4}\frac{1/\omega}{\left\lfloor 1/\omega\right\rfloor
}\right)  ^{p}.
\]

\noindent\underline{\textbf{Eighth Step.}}\textbf{\ General case: estimation
of the ratio.}

Combine the last inequality with inequalities (\ref{ginum}) - (\ref{gidenom2})
to obtain%

\begin{align*}
r  &  \geq\frac{\left(  1-\epsilon\right)  \dfrac{\sin\pi\omega}{\pi}%
q^{1-1/p}\delta_{p}^{1/p}\left(  mQ\right)  ^{1-1/p}}{Cq^{2-2/p}+\Delta
q^{1-1/p}\omega2^{1+1/p}\left(  \left\lfloor 1/\omega\right\rfloor
+1+\dfrac{\left\lfloor 1/\omega\right\rfloor }{p-1}\rho\right)  ^{1/p}}\\
&  =\frac{\left(  1-\epsilon\right)  \dfrac{\sin\pi\omega}{\pi}}{\frac
{C}{\delta_{p}^{1/p}}\left(  \frac{q}{mQ}\right)  ^{1-1/p}+\omega
2^{1+1/p}\left(  \left\lfloor 1/\omega\right\rfloor +1+\dfrac{\left\lfloor
1/\omega\right\rfloor }{p-1}\rho\right)  ^{1/p}+o\left(  1\right)  }\\
\  &  \geq\frac{\left(  1-\epsilon\right)  \dfrac{\sin\pi\omega}{\pi}%
}{\epsilon+\omega2^{1+1/p}\left(  \left\lfloor 1/\omega\right\rfloor
+1+\frac{\left\lfloor 1/\omega\right\rfloor }{p-1}\rho\right)  ^{1/p}+o\left(
1\right)  },
\end{align*}
where we have used $\Delta\simeq\delta_{p}^{1/p}(mQ)^{1-1/p}$. Here finally is
the choice of $Q$: since $q/m\leq1$, $Q=Q\left(  \epsilon,p\right)  $ is
chosen to make the first term of the denominator less than $\epsilon$. Now
since $\epsilon$ was arbitrary, we can take it to be zero and then we can take
$q$ as large as we need to get rid of the $o\left(  1\right)  $ term. In other
words,we have%

\begin{equation}
C_{p}\geq\frac{\sin\pi\omega/\left(  \pi\omega\right)  }{2^{1+1/p}\left(
\left\lfloor 1/\omega\right\rfloor +1+\frac{\left\lfloor 1/\omega\right\rfloor
}{p-1}\left(  \frac{1}{4}\frac{1/\omega}{\left\lfloor 1/\omega\right\rfloor
}\right)  ^{p}\right)  ^{1/p}}. \label{giratio1}%
\end{equation}
Since the numbers $\omega$ in the right hand side of inequality
(\ref{giratio1}) are easily seen to be everywhere dense in $\left[
0,1/2\right]  $, this ends the proof of Part I of Theorem 1. \smallskip

A little more was proved than what was stated in terms of the constant $c_{p}
$, in fact:

\begin{remark}
It is possible to get numerical estimates for particular values of $p$ by
picking a value of $\omega$ that maximizes the right hand side of this
inequality. For example, if $p=2,$ then setting $\omega=.34$ produces
$c_{2}=.13,$ which compares reasonably well with the known fact that
$C_{2}=.48....$ Since $\omega\in(0,.5),$ $1/\omega\in(2,\infty)$ so that
$\rho\leq(3/8)^{p}<1.$ (In the statement of the theorem we have very slightly
degraded estimate (\ref{giratio1}) by substituting $(3/4)^{p}$ for $\rho.$)
Hence as $p\rightarrow\infty,$ the denominator of the right hand side tends to
$1/2$ and therefore the right side becomes
\[
\frac{1}{2}\frac{\sin\pi\omega}{\pi\omega}%
\]
which may be made as close to $1$ as you like by picking $\omega$ small
enough. In other words, as $p$ tends to $\infty$, $c_{p}$ tends to $1/2.$ Also
our estimate, if sharp, would show that the constant is $O(p-1)$ as
$p\searrow1$, which would be consistent with our conjecture that concentration
fails for $L^{1}.$
\end{remark}

\section{Proof of Theorem 1; Part II: $C_{p}\geq c_{p}^{\ast}$ (for all
$p\geq2$)}

\noindent First pick $\omega=\omega(p)$ so that
\[
\frac{\sin^{p}\pi\omega}{\omega^{p-1}}=\sup_{t\in(0,1/2]}\frac{\sin^{p}\pi
t}{t^{p-1}}.
\]
Let $\xi$ be an irrational point of density of $E$ and fix $\epsilon>0.$ Pick
$\delta>0$ so small that $J:=[\xi-\delta,\xi+\delta](\operatorname*{mod}1)$
satisfies
\begin{equation}
|J\backslash E|<\epsilon|J|. \label{density}%
\end{equation}
Set $S_{\theta}:=\{n$ integer: $1\leq n\leq N,\left\Vert n\xi-\theta
\right\Vert \leq\frac{\omega}{2}\}$ and $f_{\theta}:=\sum_{n\in S_{\theta}%
}e(nx).$ Equation (7) on page 901 of \cite{AJS} asserts that
\begin{equation}
\int_{0}^{1}|f_{\theta}(x)|^{2}d\theta\geq\left(  \frac{\sin\pi\omega}{\pi
}\right)  ^{2}|D_{N}(x-\xi)|^{2}. \label{L2}%
\end{equation}
Combining this with H\"{o}lder's inequality,
\[
\left(  \int_{0}^{1}|f_{\theta}(x)|^{p}d\theta\right)  ^{1/p}\left(  \int
_{0}^{1}1^{p^{\prime}}d\theta\right)  ^{1/p^{\prime}}\geq\left(  \int_{0}%
^{1}|f_{\theta}(x)|^{2}d\theta\right)  ^{1/2},
\]
we get
\[
\int_{0}^{1}|f_{\theta}(x)|^{p}d\theta\geq\left(  \int_{0}^{1}|f_{\theta
}(x)|^{2}d\theta\right)  ^{p/2}\geq\left(  \frac{\sin\pi\omega}{\pi}\right)
^{p}|D_{N}(x-\xi)|^{p}.
\]
Now integrate this in $x$ over $J$ to get
\begin{equation}
\int_{0}^{1}\int_{J}|f_{\theta}(x)|^{p}dx\,d\theta\geq\left(  \frac{\sin
\pi\omega}{\pi}\right)  ^{p}\int_{\xi-\delta}^{\xi+\delta}|D_{N}(x-\xi
)|^{p}dx=\left(  \frac{\sin\pi\omega}{\pi}\right)  ^{p}\int_{-\delta}^{\delta
}|D_{N}(u)|^{p}du. \label{Lp}%
\end{equation}

Recall from the lemma of Section 2 that%
\[
\int_{0}^{1}|D_{N}(u)|^{p}du=\ell_{p}N^{p-1}+R_{p}\left(  N\right)  ,
\]
where $\ell_{p}=\left(  2/\pi\right)  \int_{0}^{\infty}\left\vert \sin
x/x\right\vert ^{p}dx$ and%
\[
R_{p}\left(  N\right)  =\left\{
\begin{array}
[c]{lll}%
O_{p}\left(  N^{p-3}\right)  &  & \text{ if }p>3\\
O\left(  \log N\right)  &  & \text{ if }p=3\\
O_{p}\left(  1\right)  &  & \text{ if }1<p<3
\end{array}
.\right.
\]

We can now make the estimate
\begin{align*}
\int_{-\delta}^{\delta}|D_{N}(u)|^{p}du  &  =\int_{-1/2}^{1/2}-2\int_{\delta
}^{1/2}\\
&  =\ell_{p}N^{p-1}-\frac{2}{\pi^{p}}\int_{\delta}^{\pi/2}\left\vert \sin N\pi
u\right\vert ^{p}/u^{p}du+R_{p}\left(  N\right) \\
&  \geq\ell_{p}N^{p-1}-\frac{2}{\pi^{p}}\int_{\delta}^{\pi/2}1/u^{p}%
du+R_{p}\left(  N\right) \\
&  =\ell_{p}N^{p-1}-\frac{2}{(p-1)\pi^{p}}\left(  1/\delta\right)
^{p-1}+R_{p}\left(  N\right)  .
\end{align*}

Putting this into estimate (\ref{Lp}) yields
\[
\int_{0}^{1}\int_{J}|f_{\theta}(x)|^{p}dxd\theta\geq\left(  \frac{\sin
\pi\omega}{\pi}\right)  ^{p}\left(  \ell_{p}N^{p-1}-\frac{2}{(p-1)\pi^{p}%
}\left(  1/\delta\right)  ^{p-1}+R_{p}\left(  N\right)  \right)  .
\]
Hence there must be at least one $\theta$ for which
\begin{equation}
\int_{J}|f_{\theta}(x)|^{p}dx\geq\left(  \frac{\sin\pi\omega}{\pi}\right)
^{p}\left(  \ell_{p}N^{p-1}-\frac{2}{(p-1)\pi^{p}}\left(  \frac{1}{\delta
}\right)  ^{p-1}+R_{p}\left(  N\right)  \right)  . \label{bigpnum}%
\end{equation}
Next observe that
\begin{equation}
\operatorname*{card}S_{\theta}=N\omega+\epsilon(N)N, \label{card}%
\end{equation}
where $\epsilon(N)\rightarrow0$ as $N\rightarrow\infty.$ (To see this one can,
for example, trace through the proof of Weyl's theorem given on pages 11--13
of K\"{o}rner's \textit{Fourier Analysis }\cite{Ko}. When the interval $[2\pi
a,2\pi b]$ appearing there is translated, the functions $f_{+}$ and $f_{-}$
are also. But translating a function amounts to multiplying its Fourier
coefficients by factors of modulus 1, whence it is easy to see that all of the
estimates depend only on $b-a$ and not on the value of $a$.) It follows that
\[
\int_{0}^{1}|f_{\theta}(x)|^{p}dx=\int_{0}^{1}|f_{\theta}(x)|^{p-2}|f_{\theta
}(x)|^{2}dx
\]%
\[
\leq\left(  N\omega+\epsilon(N)N\right)  ^{p-2}\int_{0}^{1}|f_{\theta}%
(x)|^{2}dx
\]%
\[
=\left(  N\omega+\epsilon(N)N\right)  ^{p-2}\left(  N\omega+\epsilon
(N)N\right)  =\left(  N\omega\right)  ^{p-1}+\epsilon_{1}(N)N^{p-1},
\]
where $\epsilon_{1}(N)\rightarrow0$ as $N\rightarrow\infty.$ Thus
\begin{equation}
\int_{0}^{1}|f_{\theta}(x)|^{p}dx\leq\left(  N\omega\right)  ^{p-1}%
+\epsilon_{1}(N)N^{p-1}. \label{bigpdenom}%
\end{equation}
It also follows from relations (\ref{density}) and (\ref{card}) that
\[
\left.  \int_{E}|f_{\theta}(x)|^{p}dx\right/  \int_{0}^{1}|f_{\theta}%
(x)|^{p}dx\geq\left.  \int_{J\cap E}|f_{\theta}(x)|^{p}dx\right/  \int_{0}%
^{1}|f_{\theta}(x)|^{p}dx
\]%
\[
=\left.  \int_{J}|f_{\theta}(x)|^{p}dx\right/  \int_{0}^{1}|f_{\theta}%
(x)|^{p}dx-\left.  \int_{J\backslash E}|f_{\theta}(x)|^{p}dx\right/  \int
_{0}^{1}|f_{\theta}(x)|^{p}dx.
\]

Denote the last two ratios by $I$ and $II$ respectively. We complete the proof
by showing that $I$ is big and that $II$ is small. To estimate $I$ we use
relations (\ref{bigpnum}) and (\ref{bigpdenom}).
\[
I\geq\frac{\left(  \left(  \sin\pi\omega\right)  /\pi\right)  ^{p}\left(
\ell_{p}N^{p-1}-2\left(  p-1\right)  ^{-1}\pi^{-p}\left(  1/\delta\right)
^{p-1}+R_{p}\left(  N\right)  \right)  }{\left(  N\omega\right)
^{p-1}+\epsilon_{1}(N)N^{p-1}}
\]%
\[
=\frac{\sin^{p}\pi\omega}{\pi^{p}\omega^{p-1}}\left(  \ell_{p}-2\left(
p-1\right)  ^{-1}\pi^{-p}\left(  \frac{1}{N\delta}\right)  ^{p-1}\right)
+\epsilon_{2}(N),
\]
where $\epsilon_{2}(N)\rightarrow0$ as $N\rightarrow\infty.$ Since
$[0,1]\supset J,$ we may use the estimate (\ref{bigpnum}) for the denominator
of $II,$ obtaining
\[
II<\frac{\epsilon|J|\sup_{x}|f_{\theta}(x)|^{p}}{\left(  \left(  \sin\pi
\omega\right)  /\pi\right)  ^{p}N^{p-1}\left(  \ell_{p}-2\left(  p-1\right)
^{-1}\pi^{-p}\left(  1/\left(  N\delta\right)  \right)  ^{p-1}+o(1)\right)  }
\]
so using $|J|=2\delta$ and the estimate (\ref{card}), we obtain
\[
II<\frac{2\epsilon N\delta\omega}{\left(  \sin^{p}\pi\omega\right)  \pi
^{-p}\omega^{1-p}\left(  \ell_{p}-2(p-1)^{-1}\pi^{-p}\left(  1/\left(
N\delta\right)  \right)  ^{p-1}\right)  }+\epsilon_{3}(N),
\]
where $\epsilon_{3}(N)\rightarrow0$ as $N\rightarrow\infty.$ Combine the
estimates for $I$ and $II$ to achieve
\[
C_{p}^{\ast}\geq\frac{\sin^{p}\pi\omega}{\pi^{p}\omega^{p-1}}\left(  \ell
_{p}-2(p-1)^{-1}\pi^{-p}\left(  1/\left(  N\delta\right)  \right)
^{p-1}\right)
\]%
\[
-\frac{2\epsilon N\delta\omega}{\left(  \sin^{p}\pi\omega\right)  \pi
^{-p}\omega^{1-p}\left(  \ell_{p}-2(p-1)^{-1}\pi^{-p}\left(  1/\left(
N\delta\right)  \right)  ^{p-1}\right)  }-\epsilon_{3}(N).
\]

Given any $\eta>0,$ find $M$ so large that $2(p-1)^{-1}\pi^{-p}\left(
1/\left(  N\delta\right)  \right)  ^{p-1}<\eta,$ whenever $N\delta>M.$ Then
pick $\epsilon>0$ so small that
\[
\frac{2\epsilon(M+1)\omega}{\left(  \sin^{p}\pi\omega\right)  \pi^{-p}%
\omega^{1-p}\left(  \ell_{p}-\eta\right)  }<\eta.
\]
Next pick $\delta>0$ so small that estimate (\ref{density}) holds for this
$\epsilon.$ Finally choose $N$ so that $M<N\delta<M+1.$ It then follows from
our last estimate for $C_{p}^{\ast}$ that
\[
\left(  C_{p}^{\ast}\right)  ^{p}\geq\frac{\sin^{p}\pi\omega}{\pi^{p}%
\omega^{p-1}}\left(  \ell_{p}-\eta\right)  -\eta-\epsilon_{3}(N).
\]
Since $\eta\,$ was arbitrary and since increasing $M$ also shrinks
$\epsilon_{3}(N)$,
\[
\left(  C_{p}^{\ast}\right)  ^{p}\geq\frac{\sin^{p}\pi\omega}{\pi^{p}%
\omega^{p-1}}\ell_{p}=\frac{\sin^{p}\pi\omega}{\pi^{p}\omega^{p-1}}\frac
{2}{\pi}\int_{0}^{\infty}|\frac{\sin x}{x}|^{p}dx.
\]
Thus Part II of Theorem 1 is proved.

\section{Does concentration fail when $p=1$?\medskip}

\noindent\textbf{Conjecture 1. }\textit{Concentration fails for} $L^{1}.$
\textit{More specifically, there is an absolute constant} $D$ \textit{such
that if} $J=[\frac{1}{q}-\frac{1}{mq},\frac{1}{q}+\frac{1}{mq}],$
\textit{where} $m>q^{2},$ \textit{then for every idempotent} $f,$%
\begin{equation}
\int_{J}|f|dx\left/  \int_{0}^{1}|f|dx\right.  \leq\frac{D}{\ln q}. \label{L1}%
\end{equation}

Define a \textit{special idempotent }to be a idempotent of the form
$D_{k}(lx)D_{m}(x)$ for positive integers $k,l,$and $m.$ Our main reason for
believing this conjecture is that (1) estimates of the type made in section 2
above show that the conjecture holds for all special idempotents and (2) the
special idempotents do provide both the correct asymptotic behavior at
$p=\infty$ and also the exact maximizing constant when $p=2$.\medskip

To see that this last point is so, we must sharpen the estimates that we made
in section 2 above. Set $p=2$ and use the exact calculation
\[
\int_{0}^{1}|D_{mq}(qx)D_{\omega q}(x)|^{2}dx=m\omega q^{2}
\]
for the denominator in the quantity (\ref{ratio}), estimate the numerator as
was done in section 2, and then replace $d_{mq}$ by its exact value $mq;$ then
our estimate for (\ref{ratio}) is improved to
\[
\frac{mq\left(  \left(  \sin\pi\omega\right)  /\pi\right)  ^{2}q}%
{mq\cdot\omega q}=\frac{\sin{}^{2}\pi\omega}{\pi\cdot\pi\omega},
\]
which shows that, for the appropriate choice of $\omega,$ the best possible
constant is achievable even if the supremum is taken only over the small
subclass of special idempotents.\medskip

Conjecture 1 is supported even more strongly by evidence that the following
conjecture might be true.

\noindent\textbf{Conjecture 2. }\textit{Let} $\wp_{n}:=\{\sum_{k=0}%
^{n-1}\epsilon_{k}e(kx):\epsilon_{k}$ \textit{is} $0$ \textit{or} $1\}.$
\textit{Then there is an absolute constant }$c$\textit{\ such that for every}
$n,$%
\[
\sup_{f\in\wp_{n}}\int_{\frac{1}{q}-\frac{1}{mq}}^{\frac{1}{q}+\frac{1}{mq}%
}|f(x)|dx\left/  \int_{0}^{1}|f(x)|dx\right.  \leq c\cdot\sup_{%
\begin{array}
[c]{c}%
f\in\wp_{n}\\
f\text{ special}%
\end{array}
}\int_{\frac{1}{q}-\frac{1}{mq}}^{\frac{1}{q}+\frac{1}{mq}}|f(x)|dx\left/
\int_{0}^{1}|f(x)|dx\right.  .
\]

We should remark that this conjecture is trivial when $n\leq2,$ since the
smallest idempotent that is not special is $e(0x)+e(1x)+e(3x).$ It is easy to
see that conjecture 2 easily implies conjecture 1. Supporting numerical
evidence for conjecture 2 consists primarily of the fact that for those values
of $q$ we have looked at, the vast majority of the \textquotedblright
best\textquotedblright\ functions (i.e. functions which produced the largest
ratio) were special. (If this were always the case, conjecture 2 would hold
with $c=1,$ which would be a spectacular result.) However, there are some
non-special functions which do beat out the special functions for certain
values of $n$, hence the need for the constant $c$. For instance, when $q=6$
the special functions were not always found to be the best. In particular, the
function
\[
f=e(0x)+e(1x)+e(5x)+e(6x)+e(7x)+e(12x)
\]
produces the largest ratio for $n\leq13,$ yet this is not a special function.
But the best special function for $n\leq13$ is
\[
D_{2}(6x)D_{3}(x)=e(0x)+e(1x)+e(2x)+e(6x)+e(7x)+e(8x),
\]
whose ratio is only $.98$ of the ratio produced by $f$. Table 1 lays out the
smallest values of $c$ observed for various values of $q$.\newpage
\ \begin{table}[h]
\centering
\begin{tabular}
[c]{|l|l|c|c|c|}\hline
\textbf{q} & \textbf{c} & \textbf{n's studied} & \textbf{minimum c at} &
\textbf{best projection at highest n}\\\hline
2 & 1 & 3,\dots,14 & -- & $D_{8}(2x)$\\
3 & 1 & 3,\dots, 15 & -- & $D_{6}(3x)$\\
4 & 1 & 3,\dots,16 & -- & $D_{2}(x)D_{4}(4x)$\\
5 & .87 & 3,\dots,16 & n=15 & $D_{7}(x)(1-e(2x)+e(4x))$\\
6 & .98 & 3,\dots,18 & n=13 & $D_{3}(x)D_{3}(6x)$\\
7 & .87 & 3,\dots,16 & n=15 & $D_{2}(x)D_{3}(7x)$\\
8 & .97 & 3,\dots,22 & n=18 & $D_{4}(x)D_{3}(8x)$\\
10 & .98 & 3,\dots,22 & n=22 & $D_{3}(x)(1+e(9x)+e(19x))+e(12x)$\\
12 & 1 & 3,\dots,20 & -- & $D_{6}(x)D_{2}(12x)$\\
16 & 1 & 3,\dots,20 & -- & $D_{6}(x)D_{2}(15x)$\\
32 & 1 & 3,\dots,16 & -- & $D_{15}(x)$\\\hline
\end{tabular}
\caption{Some observed upper bounds for c}%
\label{key}%
\end{table}\ 

\noindent It should be noted that some values of $q$ were studied to larger
values of $n$ than others. Clearly, the computation time is exponential in $n
$, so going up to say $n=22$ amounts to computing $2^{22}$ ratios of
integrals.\medskip

\textbf{Acknowledgment.\textit{\ }}\textit{We thank the referee for many
improvements. In particular, he suggested Lemma \ref{l:4} and explained to us
how to use it to guarantee concentration on sets of positive measure, rather
than just on intervals, when }$\mathit{1<p<2}$\textit{.}

B. Anderson

130 Channing Ln

Chapel Hill, NC 27516

USA

\textit{E-mail address}: bruce.b.anderson@csfb.com\medskip

J. M. Ash

Department of Mathematical Sciences

DePaul University

Chicago, IL 60614

USA

\textit{E-mail address}: mash@math.depaul.edu\medskip

R. Jones

Conserve School

5400 N. Black Oak Lake Drive

Land O'Lakes, WI 54540 (USA)

\textit{E-mail address}: Roger.Jones@ConserveSchool.org\medskip

D. G. Rider

Department of Mathematics

University of Wisconsin

480 Lincoln Drive

Madison, WI 53706-1313

USA

\textit{E-mail address}: drider@math.wisc.edu\medskip

B. Saffari

D\'{e}partement de Math\'{e}matiques

Universit\'{e} de Paris XI (Orsay)

91405 Orsay Cedex

France

\textit{E-mail address}: Bahman.Saffari@math.u-psud.fr
\end{document}